\documentclass[preprint,12pt]{elsarticle}

\usepackage{amssymb}
\usepackage{lineno}

\biboptions{sort&compress}

\journal{For my students.}

\usepackage{amsmath}
\DeclareMathAlphabet{\mathbfit}{\encodingdefault}{\rmdefault}{bx}{sl}
\DeclareMathAlphabet{\mathslsf}{\encodingdefault}{\sfdefault}{m}{sl}
\DeclareMathAlphabet{\mathbfsf}{\encodingdefault}{\sfdefault}{bx}{n}
\usepackage{algorithm}
\usepackage{algpseudocode}
\usepackage{longtable}
\newcommand{\textfrac}[2]
   {\mbox{\leavevmode \kern.1em\raise.5ex\hbox{\ensuremath{\scriptstyle{#1}}}\kern-.1em\hspace{0em}/\hspace{0em}\kern-.1em\lower.25ex\hbox{\ensuremath{\scriptstyle{#2}}}\kern.1em
   }}
\newcommand{\scriptfrac}[2] 
   {\mbox{\leavevmode \kern.05em\raise.25ex\hbox{\ensuremath{\scriptscriptstyle{#1}}}\kern-.05em\hspace{0em}\hbox{\ensuremath{\scriptstyle{/}}}\hspace{0em}\kern-.05em\lower.125ex\hbox{\ensuremath{\scriptscriptstyle{#2}}}\kern.05em
   }}
\usepackage{mathabx}
\newsavebox{\accentbox}

\newcommand{\ibar}{\mbox{
   \lower.35ex\hbox{\ensuremath{\mathchar'26}}\kern-1ex\raise.35ex\hbox{\ensuremath{\int}}\kern0ex
   }}
\newcommand{\tint}{\hbox{\ensuremath{\int}}\kern0ex}
\begin{document}

\begin{frontmatter}



\title{A Technical Note: \\
	 Two-Step PECE Methods for Approximating Solutions \\
	 To First- and Second-Order ODEs}


\author[add]{A.\ D.~Freed}%\corref{corr}}

\ead{afreed@tamu.edu}

%\cortext[corr]{Corresponding author}

\date{\today}

\address[add]{Department of Mechanical Engineering, 
              Texas A\&\mbox{M} University,
              College Station, \\ TX 77843, 
              United States}
              
\begin{abstract}

Two-step predictor\slash corrector methods are provided to solve three classes of problems that present themselves as systems of ordinary differential equations (ODEs).  In the first class, velocities are given from which displacements are to be solved.  In the second class, velocities and accelerations are given from which displacements are to be solved.  And in the third class, accelerations are given from which velocities and displacements are to be solved.  Two-step methods are not self starting, so compatible one-step methods are provided to take that first step with.  An algorithm is presented for controlling the step size so that the local truncation error does not exceed a specified tolerance.

\end{abstract}

%\begin{keyword}

% \sep \sep \sep \sep \sep 


% \MSC 15A06 \sep 65F05 \sep 74A05 \sep 74A10 \sep 74A15 \sep 74A20

%\end{keyword}

\end{frontmatter}

% \linenumbers


\section{Multi-Step Methods for Mechanical Engineers}

Multi-step methods \cite{Butcher08,HairerWanner91} are numerical schemes that are used to approximate solutions for systems of ODEs which commonly arise in engineering practice.  Because the intended readers of this document are my students, whom will become Mechanical Engineers upon graduation, I present these methods using variables that are intuitive to them: time $t$ is the independent variable of integration, and position $\mathbf{x} = \{ x_1 , x_2 , x_3 \}^{\mathsf{T}}$ is the dependent variable of integration (plus, sometimes, velocity) while velocity $\mathbf{v} = \{ v_1 , v_2 , v_3 \}^{\mathsf{T}} = \mathbf{v}(t, \mathbf{x})$ and acceleration $\mathbf{a} = \{ a_1 , a_2 , a_3 \}^{\mathsf{T}} = \mathbf{a}(t, \mathbf{x}, \mathbf{v})$ are functions of these independent and dependent variables.  The time rate-of-change of acceleration is jerk $\dot{\mathbf{a}} = \{ \dot{a}_1 , \dot{a}_2 , \dot{a}_3 \}^{\mathsf{T}}$, which is introduced as a means by which improvements in solution accuracy can be made.  Problems like these commonly arise in applications within disciplines like kinematics, dynamics, thermo\-dynamics, vibrations, controls, process kinetics, etc.  The methods presented in this document apply to systems of any dimension, it is just that $t$, $\mathbf{x}$, $\mathbf{v}$ and $\mathbf{a}$ are physical notions for which my students have intuitive understanding.

Current engineering curricula expose students to some basic methods like Euler's method (you should never use forward-Euler by itself), a simple Euler predictor with a trapezoidal corrector, often called Heun's method, and \textit{the\/} Runge-Kutta method. Kutta \cite{Kutta01} derived \textit{the\/} Runge-Kutta method (Runge played no part here).  This method, likely the most popular of all ODE solvers, was not the method Kutta actually advocated for use.  He derived a more accurate fourth-order method in his paper---a method that has sadly become lost to the obscurity of dusty shelves.  

The intent of this note is to inform my students about the existence and utility of a whole other class of ODE solvers that have great value in many applications.  These are called multi-step methods.  They make an informed decision on the direction that its solution will advance into the future based upon where it has been in the recent past.  In contrast, Runge-Kutta methods sample multiple paths in the present to make an informed decision on the direction that its solution will advance into the future.  The past does not enter into the Runge-Kutta process.  These two classes of numerical methods are fundamentally different in this regard.  There is an emerging field within computational mathematics where these two approaches are being melded into one.  They are called general linear methods, and two such methods can be found in Appendix~D of my textbook \cite{Freed14}.  We will not address them here.

\section{The Objective}

Throughout this document we shall consider an interval in time $[0,T]$ over which $N$ solutions are to be extracted at nodes $n = 1,2, \ldots, N$ spaced at uniform intervals in time with a common step size of $h = T/N$ separating them.  This is referred to as the global step size.  A local step size will be introduced later, which will be the actual step size that an integrator uses to advance along its solution path.  This size dynamically adjusts to maintain solution accuracy, and is under the control of a proportional integral (PI) controller.  

Node $n$ is located at current time.  Here is where the solution front resides.  Node $n \! - \! 1$ is where the previous solution was acquired, while node $n \! + \! 1$ is the where the next solution is to be calculated.  In this regard,  information storage required by these methods is compatible with memory strategies and coding practices adopted by many industrial codes like finite elements.  This requirement of working solely with nodes $n \! - \! 1$, $n$, $n \! + \! 1$ will limit the accuracy that one can achieve with these methods.  Higher-order multi-step methods require more nodes, and as such, more information history.

Our objective is to construct a collection of numerical methods that resemble the popular, second-order, backward-difference formula \cite{HairerWanner91} denoted as BDF2 in the literature and software packages.  BDF2 is described by
\begin{equation*}
\mathbf{x}_{n+1} = \tfrac{1}{3} \bigl( 4 \mathbf{x}_n - 
\mathbf{x}_{n-1} \bigr) + \tfrac{2}{3} \, h
\mathbf{v}_{n+1} + \mathcal{O} (h^3) 
\end{equation*}
and is an implicit method in that $\mathbf{v} = \mathbf{v} (t, \mathbf{x})$, typically, and therefore $\mathbf{x}_{n+1}$ appears on both sides of the equals sign.  There are good reasons for selecting this numerical model upon which to construct other methods; specifically, BDF2 is a convergent method in that it is consistent and A~stable \cite{Butcher08}.  These are noble properties to aspire to, but whose discussion lies beyond the scope of this document.  

Here your professor seeks to provide techniques that address three questions: \textit{i\/}) How can one apply an implicit multi-step method where you need to know the solution to get the solution? \textit{ii\/}) How can one startup a multi-step method, because at the initial condition there is no solution history? and \textit{iii\/}) Numerical ODE solvers typically solve first-order systems, but Newton's Laws for Motion are described with a second-order system.  How can one construct an ODE solver designed to handle these types of problems?

An answer to the first question is: We will introduce a predictor to get an initial solution estimate; specifically, predict\slash evaluate\slash correct\slash evalute (PECE) schemes are developed.  An answer to the second question is: A single-step method can be used to start up a two-step method.  And an answer to the third question is: We will use the natural features of multi-step methods and Taylor series expansions to construct solvers for second-order ODEs.  Several of the methods found in this document are not found in the literature.  Your professor created them just for you!

\subsection{Strategy}

The strategy used to construct multi-step algorithms is to expand an appropriate linear combination of Taylor series for displacement $\mathbf{x}$ taken about solution nodes at discrete times.  In our case, expansions are taken about times $t_{n-1}$, $t_n$ and $t_{n+1}$ such that their sum replicates the general structure of the BDF2 method.  Specifically, we seek two-step methods with constituents $\mathbf{x}_{n+1} = \tfrac{1}{3} (4 \textbf{x}_n - \textbf{x}_{n-1}) + \cdots$ that are common betwixt them.  

Each Taylor series is expanded out to include acceleration $\mathbf{a}$ for methods that solve first-order ODEs, and each Taylor series is expanded out to include jerk $\dot{\mathbf{a}}$ for methods that solve second-order ODEs.  The pertinent series for displacement include
\begin{subequations}
	\label{TaylorDisplacements}
	\begin{align}
	\mathbf{x}_{n+1} & = \mathbf{x}_n + h \mathbf{v}_n + 
		\tfrac{1}{2} h^2 \mathbf{a}_n + \tfrac{1}{6} h^3 
		\dot{\mathbf{a}}_n + \cdots 
		\label{displacementA} \\
	\mathbf{x}_n & = \mathbf{x}_{n+1} - h \mathbf{v}_{n+1} + 
		\tfrac{1}{2} h^2 \mathbf{a}_{n+1} - 
		\tfrac{1}{6} h^3 \dot{\mathbf{a}}_{n+1} + \cdots 
		\label{displacementB} \\
	\mathbf{x}_n & = \mathbf{x}_{n-1} + h \mathbf{v}_{n-1} + 
		\tfrac{1}{2} h^2 \mathbf{a}_{n-1} + 
		\tfrac{1}{6} h^3 \dot{\mathbf{a}}_{n-1} + \cdots 
		\label{displacementC} \\
		\mathbf{x}_{n-1} & = \mathbf{x}_n - h \mathbf{v}_n + 
		\tfrac{1}{2} h^2 \mathbf{a}_n - \tfrac{1}{6} h^3 
		\dot{\mathbf{a}}_n + \cdots 
		\label{displacementD}
	\end{align}
\end{subequations}
where the set of admissible expansions only involve nodes $n \! - \! 1$, $n$ and $n \! + \! 1$.  Once these are in place, like Taylor expansions for the velocity are secured 
\begin{subequations}
	\label{TaylorVelocities}
	\begin{align}
	\mathbf{v}_{n+1} & = \mathbf{v}_n + h \mathbf{a}_n + 
	\tfrac{1}{2} h^2 \dot{\mathbf{a}}_n + \cdots 
	\label{velocityA} \\
	\mathbf{v}_n & = \mathbf{v}_{n+1} - h \mathbf{a}_{n+1} + 
	\tfrac{1}{2} h^2 \dot{\mathbf{a}}_{n+1} + \cdots 
	\label{velocityB} \\
	\mathbf{v}_n & = \mathbf{v}_{n-1} + h \mathbf{a}_{n-1} + 
	\tfrac{1}{2} h^2 \dot{\mathbf{a}}_{n-1} + \cdots 
	\label{velocityC} \\
	\mathbf{v}_{n-1} & = \mathbf{v}_n - h \mathbf{a}_n + 
	\tfrac{1}{2} h^2 \dot{\mathbf{a}}_n + \cdots .
	\label{velocityD}
	\end{align}
\end{subequations}
These series are solved for acceleration for the first-order ODE solvers, and for jerk for the second-order ODE solvers.  These solutions for acceleration\slash jerk are then inserted back into the original series for displacement.  The net effect is to incorporate contributions for acceleration\slash jerk by approximating them in terms of velocities and, possibly, accelerations, thereby increasing the order of accuracy for the overall method by one order, e.g. from second-order, i.e., $\mathcal{O}(h^3)$, to third-order, viz., $\mathcal{O}(h^4)$, for the second-order ODE methods.  This is accomplished without the solver explicitly needing any information about jerk from the user, which would be hard to come by in practice.

We speak of a method being, say, second-order accurate, and designate this with the notation $\mathcal{O}(h^3)$.  There may seem to be an apparent discrepancy between the order of a method and the exponent of $h$.  This comes into being because the `order' of a method represents the global order of accuracy in a solution, whereas the exponent on the $h$-term represents an order of accuracy in the solution over a local step of integration with the exponent on $h$ designating the order of its error estimate.

Our objective, viz., $\mathbf{x}_{n+1} = \mathbf{x}_n + \cdots$ for one-step (startup) methods and  $\mathbf{x}_{n+1} = \tfrac{1}{3} (4 \textbf{x}_n - \textbf{x}_{n-1}) + \cdots$ for two-step methods, is achieved by applying the following linear combinations of Taylor series
\begin{displaymath}
   \textrm{predictors} \Leftarrow \begin{cases}
      1 (\ref{displacementA}) & \textrm{one-step} \\
      1 (\ref{displacementA}) - 
      \tfrac{1}{6} (\ref{displacementC}) +
      \tfrac{1}{6} (\ref{displacementD}) & 
      \textrm{two-step}
   \end{cases}
\end{displaymath}
and
\begin{displaymath}
	\textrm{correctors} \Leftarrow \begin{cases}
	   \tfrac{1}{2} (\ref{displacementA}) - 
	   \tfrac{1}{2} (\ref{displacementB}) & \textrm{one-step} \\
	   \tfrac{4}{3} (\ref{displacementA}) +
	   \tfrac{1}{3} (\ref{displacementB}) +
	   \tfrac{1}{3} (\ref{displacementD}) & 
	   \textrm{two-step } \# 1 \\
	   \tfrac{4}{3} (\ref{displacementA}) +
	   \tfrac{1}{3} (\ref{displacementB}) -
	   \tfrac{1}{6} (\ref{displacementC}) +
	   \tfrac{1}{6} (\ref{displacementD}) & 
	   \textrm{two-step } \# 2
	\end{cases}
\end{displaymath}
with
\begin{displaymath}
	\textrm{truncation errors} \Leftarrow \begin{cases}
	\tfrac{1}{2} \| (\ref{displacementA}) +
	(\ref{displacementB}) \| & \textrm{one-step} \\
	\tfrac{1}{6} \| 2 (\ref{displacementA}) +
	2 (\ref{displacementB}) + 1 (\ref{displacementC}) +
	1 (\ref{displacementD}) \| & \textrm{two-step } \# 1 \\
	\tfrac{1}{3} \| (\ref{displacementA}) +
	(\ref{displacementB}) \| & \textrm{two-step } \# 2
\end{cases}
\end{displaymath}
wherein the parenthetical numbers refer to the sub-equations listed in Eq.~(\ref{TaylorDisplacements}) and where the coefficients out front designate the weight applied to that formula.  To be a corrector requires expansion (\ref{displacementB}), which must not appear in a predictor.

There are two ways to construct a corrector that satisfy our conjecture, and both will be used. A design objective is to come up with a predictor\slash corrector pair that weigh their contributions the same; specifically, their displacements are weighted the same, their velocities are weighted the same, and when present, their accelerations are weighted the same, too.

\section{PECE Methods for First-Order ODEs}
\label{Sec:firstOrder}

The following algorithm is suitable for numerically approximating solutions to stiff systems of ODEs, which engineers commonly encounter.  The idea of mathematical stiffness is illustrated through an example in \S\ref{Sec:Brusselator}.

For this class of problems it is assumed that velocity is described as a function in time and displacement, e.g., at step $n$ a formula would give $\mathbf{v}_n = \mathbf{v} (t_n , \mathbf{x}_n)$.  An initial condition $\mathbf{x} (0) = \mathbf{x}_0$ is required to start an analysis.  The objective is to solve this ODE for displacement $\mathbf{x}_{n+1}$ evaluated at the next moment in time $t_{n+1}$, wherein $n$ sequences as $n = 0, 1, \ldots , N \! - \! 1$.

Heun's method is used to take the first integration step.  Begin by applying a predictor (it is a forward Euler step) 
\begin{subequations}
	\label{startUp1stOrderODEs}
	\begin{align}
	   \mathbf{x}_1^p & = \mathbf{x}_0 + h \mathbf{v}_0 + 
	   \mathcal{O} (h^2)
	   \label{startUp1stOrderPredictor} \\
	   \intertext{which is to be followed with an evaluation for velocity $\mathbf{v}^p_1 = \mathbf{v} (t_1 , \mathbf{x}_1^p)$ using this predicted estimate for displacement.  A corrector is then applied (it is the trapezoidal rule)}
	   \mathbf{x}_1 & = \mathbf{x}_0 + \tfrac{1}{2} h 
	   \bigl( \mathbf{v}_1^p + \mathbf{v}_0 \bigr) + 
	   \mathcal{O} (h^3)
	   \label{startUp1stOrderCorrector}
	\end{align}
\end{subequations}
after which a final re-evaluation for velocity $\mathbf{v}_1 = \mathbf{v} (t_1 , \mathbf{x}_1)$ is made and the first step comes to a close.  In this case, using another Taylor series to subtract out the influences from acceleration did not bring about any change to the formula.  This is not unexpected, as the trapezoidal method is already second-order accurate, i.e., it has a truncation error on the order of $\mathcal{O}(h^3)$.  The step counter is assigned a value of $n = 1$, after which control of the solution process is passed over to the following method.

For entering step counts that lie within the interval $n=1$ to $n = N \! - \! 1$, numeric integration continues by employing a predictor 
\begin{subequations}
	\label{1stOrderODEs}
	\begin{align}
	\mathbf{x}_{n+1}^p & = \tfrac{1}{3} 
	\bigl( 4 \mathbf{x}_n - \mathbf{x}_{n-1} \bigr) + 
	\tfrac{2}{3} h \bigl( 2\mathbf{v}_n - \mathbf{v}_{n-1} 
	\bigr) + \mathcal{O} (h^3)
	\label{1stOrderPredictor} \\
	\intertext{followed by an evaluation for velocity via $\mathbf{v}^p_{n+1} = \mathbf{v} (t_{n+1} , \mathbf{x}_{n+1}^p)$ using this predicted estimate for displacement.  Here including correction terms for acceleration changed $\tfrac{1}{6} h ( 5 \mathbf{v}_n - \mathbf{v}_{n-1})$ to $\tfrac{2}{3} h ( 2 \mathbf{v}_n - \mathbf{v}_{n-1} )$ and in the process improved its accuracy from $\mathcal{O}(h^2)$ to $\mathcal{O}(h^3)$.  The corrector obtained according to our recipe for a type \#1 method is}
	\mathbf{x}_{n+1} & = \tfrac{1}{3} 
	\bigl( 4 \mathbf{x}_n - \mathbf{x}_{n-1} \bigr) + 
	\tfrac{2}{3} h \mathbf{v}^p_{n+1} + \mathcal{O} (h^3)
	\label{1stOrderCorrector}
	\end{align}
\end{subequations} 
which culminates with a re-evaluation for $\mathbf{v}_{n+1} = \mathbf{v} ( t_{n+1}, \mathbf{x}_{n+1})$.  This corrector is the well-known BDF2 formula, the method we are generalizing around.  Including correction terms for acceleration changed $- \tfrac{1}{3} h ( \mathbf{v}^p_{n+1} - 3 \mathbf{v}_n)$ to $\tfrac{2}{3} h \mathbf{v}^p_{n+1}$ and in the process improved its accuracy from $\mathcal{O}(h^2)$ to $\mathcal{O}(h^3)$.  For both integrators, displacement has weight 1, while velocity has weight $\tfrac{2}{3} h$.  The predictor and corrector are consistent in this regard, a required design objective when deriving an admissible PECE method.

Variables are to be updated according to $n \! - \! 1 \leftarrow n$ and $n \leftarrow n \! + \! 1$ after which counter $n$ gets incremented. After finishing with the data management, the solution is ready for advancement to the next integration step, with looping continuing until $n = N$ whereat the solution becomes complete.

\section{PECE Methods for Second-Order ODEs}
\label{Sec:secondOrder}

For this class of problems it is assumed that the velocity is described as a function of time and displacement, e.g.,  $\mathbf{v}_n = \mathbf{v} (t_n , \mathbf{x}_n)$, and likewise, the acceleration is also a prescribed function in terms of time, displacement and velocity, e.g., $\mathbf{a}_n = \mathbf{a} (t_n , \mathbf{x}_n , \mathbf{v}_n)$.  An initial condition is to be supplied by the user, viz., $\mathbf{x} (0) = \mathbf{x}_0$.  The objective of this method is to solve this second-order ODE for displacement $\mathbf{x}_{n+1}$, which is to be evaluated at the next moment in time $t_{n+1}$, wherein $n=0,1, \ldots , N \! - \! 1$.

Like the previous method, this is a two-step method so, consequently, it is not self starting.  To take a first step, apply the predictor (a straightforward Taylor series expansion) 
\begin{subequations}
	\label{startup}
	\begin{align}
	\mathbf{x}_1^p & = \mathbf{x}_0 + h \mathbf{v}_0 +
	\tfrac{1}{2} h^2 \mathbf{a}_0 + \mathcal{O} (h^3) 
	\label{startupPredictor} \\
	\intertext{followed by evaluations $\mathbf{v}^p_1 = \mathbf{v} (t_1, \mathbf{x}^p_1)$ and $\mathbf{a}^p_1 = \mathbf{a} (t_1, \mathbf{x}^p_1, \mathbf{v}^p_1)$ to prepare for executing its corrector}
	\mathbf{x}_1 & = \mathbf{x}_0 + \tfrac{1}{2} h 
	\bigl( \mathbf{v}^p_1 + \mathbf{v}_0 \bigr) -
	\tfrac{1}{12} h^2 \bigl( \mathbf{a}^p_1 - 
	\mathbf{a}_0 \bigr) + \mathcal{O} (h^4) 
	\label{startupCorrector}
	\end{align}
\end{subequations}
after which one re-evaluates $\mathbf{v}_1 = \mathbf{v} (t_1, \mathbf{x}_1)$ and $\mathbf{a}_1 = \mathbf{a} (t_1, \mathbf{x}_1, \mathbf{v}_1)$.  Including correction terms for jerk changed $- \tfrac{1}{4} h^2 ( \mathbf{a}^p_{n+1} - \mathbf{a}_n)$ to $- \tfrac{1}{12} h^2 ( \mathbf{a}^p_{n+1} - \mathbf{a}_n)$ and in the process improved its accuracy from $\mathcal{O}(h^3)$ to $\mathcal{O}(h^4)$.  After this integrator has been run once, a switch is made to employ the two-step PECE method described below to finish up.

For entering step counts that lie within the interval $n=1$ to $n = N \! - \! 1$, numeric integration continues by employing a predictor
\begin{subequations}
	\label{PECE}
	\begin{align}
	   \mathbf{x}_{n+1}^p & = \tfrac{1}{3} \bigl(
	   4 \mathbf{x}_n - \mathbf{x}_{n-1} \bigr) + 
	   \tfrac{1}{6} h \bigl( 3 \mathbf{v}_n + 
	   \mathbf{v}_{n-1} \bigr) \notag \\ 
	   \mbox{} & \hspace{4.5cm} + 
	   \tfrac{1}{36} h^2 \bigl( 31 \mathbf{a}_n - 
	   \mathbf{a}_{n-1} \bigr) + \mathcal{O} (h^4) 
	   \label{predictor} \\
	   \intertext{followed by evaluations $\mathbf{v}^p_{n+1} = \mathbf{v} (t_{n+1}, \mathbf{x}^p_{n+1})$ and $\mathbf{a}^p_{n+1} = \mathbf{a} (t_{n+1}, \mathbf{x}^p_{n+1}, \mathbf{v}^p_{n+1})$ to be made sequentially.  Here including correction terms for jerk changed $\tfrac{1}{6} h ( 5 \mathbf{v}_n - \mathbf{v}_{n-1} )$ $+$ $\tfrac{1}{12} h^2 ( 7 \mathbf{a}_n - \mathbf{a}_{n-1} )$ to $\tfrac{1}{6} h ( 3 \mathbf{v}_n + \mathbf{v}_{n-1} )$ $+$ $\tfrac{1}{36} h^2 ( 31 \mathbf{a}_n - \mathbf{a}_{n-1} )$ and in the process improved its accuracy from $\mathcal{O}(h^3)$ to $\mathcal{O}(h^4)$.  A corrector that is consistent with the above predictor is}
	   \mathbf{x}_{n+1} & = \tfrac{1}{3} \bigl(
	   4  \mathbf{x}_n - \mathbf{x}_{n-1} \bigr) +
	   \tfrac{1}{24} h \bigl( \mathbf{v}^p_{n+1} +
	   14 \mathbf{v}_n + \mathbf{v}_{n-1} \bigr) 
	   \notag \\
	   \mbox{} & \hspace{4.5cm} +
	   \tfrac{1}{72} h^2 \bigl( 10 \mathbf{a}^p_{n+1} + 
	   51 \mathbf{a}_n - \mathbf{a}_{n-1} \bigr) + 
	   \mathcal{O} (h^4)
	   \label{corrector}
	\end{align}
\end{subequations}
whose derivation follows below in \S\ref{Sec:derivation}.  With the corrector having been run, finish by re-evaluating $\mathbf{v}_{n+1} = \mathbf{v} (t_{n+1}, \mathbf{x}_{n+1})$ and $\mathbf{a}_{n+1} = \mathbf{a} (t_{n+1}, \mathbf{x}_{n+1}, \mathbf{v}_{n+1})$.

Variables are to be updated according to $n \! - \! 1 \leftarrow n$, $n \leftarrow n \! + \! 1$, plus the counter $n$ gets incremented. After that the solution is ready for advancement to the next integration step, with looping continuing until $n = N$ whereat the solution becomes complete.

\subsection{Derivation of the Corrector}
\label{Sec:derivation}

The corrector obtained via our recipe for a type \#1 corrector is
\begin{displaymath}
	\mathbf{x}_{n+1} = \tfrac{1}{3} \bigl(
	4  \mathbf{x}_n - \mathbf{x}_{n-1} \bigr) +
	\tfrac{1}{9} h \bigl( \mathbf{v}^p_{n+1} +
	5 \mathbf{v}_n \bigr) +
	\tfrac{2}{9} h^2 \bigl( \mathbf{a}^p_{n+1} + 
	3 \mathbf{a}_n \bigr) + \mathcal{O} (h^4) 
\end{displaymath} 
where inclusion of correction terms for jerk changed $\tfrac{1}{3} h ( -\mathbf{v}^p_{n+1} + 3 \mathbf{v}_n )$ $+$ $\tfrac{1}{6} h^2 ( \mathbf{a}^p_{n+1} + 5 \mathbf{a}_n )$ to $\tfrac{1}{9} h ( \mathbf{v}^p_{n+1} + 5 \mathbf{v}_n )$ $+$ $\tfrac{2}{9} h^2 ( \mathbf{a}^p_{n+1} + 3 \mathbf{a}_n )$ and in the process improved its accuracy from $\mathcal{O}(h^3)$ to $\mathcal{O}(h^4)$. 

The corrector obtained via our recipe for a type \#2 corrector is
\begin{displaymath}
    \begin{aligned}
	\mathbf{x}_{n+1} & = \tfrac{1}{3} \bigl(
	4  \mathbf{x}_n - \mathbf{x}_{n-1} \bigr) +
	\tfrac{1}{36} h \bigl( - \mathbf{v}^p_{n+1} +
	22 \mathbf{v}_n + 3 \mathbf{v}_{n-1} \bigr) \\
	\mbox{} & \hspace{4.5cm} +
	\tfrac{1}{36} h^2 \bigl( 2 \mathbf{a}^p_{n+1} + 
	27 \mathbf{a}_n - \mathbf{a}_{n-1} \bigr) + 
	\mathcal{O} (h^4) 
	\end{aligned}
\end{displaymath}
where the correction terms for jerk changed $-\tfrac{1}{6} h ( -2 \mathbf{v}_{n+1}^p + 7 \mathbf{v}_n - \mathbf{v}_{n-1} )$ $+$ $\tfrac{1}{12} h^2 ( 2 \mathbf{a}^p_{n+1} + 9 \mathbf{a}_n - \mathbf{a}_{n-1} )$ to $\tfrac{1}{36} h ( - \mathbf{v}^p_{n+1} + 22 \mathbf{v}_n + 3 \mathbf{v}_{n-1} )$ $+$ $\tfrac{1}{36} h^2 ( 2 \mathbf{a}^p_{n+1} + 27 \mathbf{a}_n - \mathbf{a}_{n-1} )$ and in the process improved its accuracy from $\mathcal{O}(h^3)$ to $\mathcal{O}(h^4)$.  

Unfortunately, neither of these two correctors is consistent with the predictor in Eq.~(\ref{predictor}).  This predictor has a weight imposed on displacement of 1, a weight imposed on velocity of $\tfrac{2}{3} h$, and a weight imposed on acceleration of $\tfrac{5}{6} h^2$.  It is desirable to seek a corrector with these same weights.  This would imply that if a field, say acceleration, were uniform over a time interval, say $[t_{n-1}, t_{n+1}]$, then both the predictor and corrector would produce the same numeric value for acceleration's contribution to the overall result at this location in time.  The correctors derived from types~\#1 and \#2 are consistent with this predictor for all contributions except acceleration.  In terms of acceleration, the predictor has a weight of $\tfrac{5}{6} h^2$, while corrector~\#1 has a weight of $\tfrac{8}{9} h^2$ and corrector \#2 has a weight of $\tfrac{7}{9} h^2$.  Curiously, averaging correctors~\#1 and \#2 does produce the correct weight. There is consistency between the predictor and this `averaged' corrector, which is the corrector put forward in Eq.~(\ref{corrector}).

\subsection{When Only Acceleration is Controlled}
\label{Sec:Newton}

There is an important class of problems that is similar to the above class in that acceleration is described through a function of state; however, velocity is not.  Velocity, like displacement, is a response function for this class of problems.  Acceleration is still described by a function of time, displacement and velocity, e.g., $\mathbf{a}_n = \mathbf{a} (t_n , \mathbf{x}_n , \mathbf{v}_n)$; however, instead of the velocity being given as a function, it, like displacement, is to be solved through integration.  Two initial conditions must be supplied, viz., $\mathbf{x} (0) = \mathbf{x}_0$ and $\mathbf{v} (0 , \mathbf{x}_0) = \mathbf{v}_0$.  This is how Newton's Second Law usually presents itself for analysis. Beeman \cite{Beeman76} constructed a different set of multi-step methods that can also be used to get solutions for this class of problems.

This is a two-step method.  Therefore, it will require a one-step method to startup an analysis.  To start integration, take the first step using predictors
\begin{subequations}
	\label{pairedStartUp}
	\begin{align}
	\mathbf{x}_1^p & = \mathbf{x}_0 + h \mathbf{v}_0 +
	\tfrac{1}{2} h^2 \mathbf{a}_0 + \mathcal{O} (h^3) 
	\label{startupDisplacementPredictor} \\
	\mathbf{v}^p_1 & = \mathbf{v}_0 + h \mathbf{a}_0 + 
	\mathcal{O} (h^3) 
	\label{startUpVelocityPredictor} \\
	\intertext{followed by an evaluation for $\mathbf{a}^p_1 = \mathbf{a} (t_1, \mathbf{x}^p_1, \mathbf{v}^p_1)$.  Their paired correctors are}
	\mathbf{x}_1 & = \mathbf{x}_0 + \tfrac{1}{2} h 
	\bigl( \mathbf{v}^p_1 + \mathbf{v}_0 \bigr) -
	\tfrac{1}{12} h^2 \bigl( \mathbf{a}^p_1 - 
	\mathbf{a}_0 \bigr) + \mathcal{O} (h^4) 
	\label{startupDisplacementCorrector} \\
	\mathbf{v}_1 & = \mathbf{v}_0 + \tfrac{1}{2} h 
	\bigl( \mathbf{a}_1^p + \mathbf{a}_0 \bigr) + 
	\mathcal{O} (h^4)
	\label{startUpVelocityCorrector}
	\end{align}
\end{subequations}
followed with a re-evaluation for $\mathbf{a}_1 = \mathbf{a} (t_1, \mathbf{x}_1, \mathbf{v}_1)$.  With the first step of integration taken, one can switch to the PECE algorithm described below.

For entering step counts that lie within the interval $n=1$ to $n = N \! - \! 1$, numeric integration continues by employing predictors
\begin{subequations}
	\label{pairedMethods}
	\begin{align}
	\mathbf{x}_{n+1}^p & = \tfrac{1}{3} \bigl(
	4 \mathbf{x}_n - \mathbf{x}_{n-1} \bigr) + 
	\tfrac{1}{6} h \bigl( 3 \mathbf{v}_n + 
	\mathbf{v}_{n-1} \bigr) \notag \\ 
	\mbox{} & \hspace{3.175cm} + 
	\tfrac{1}{36} h^2 \bigl( 31 \mathbf{a}_n - 
	\mathbf{a}_{n-1} \bigr) + \mathcal{O} (h^4) 
	\label{displacementPredictor} \\
	\mathbf{v}_{n+1}^p & = \tfrac{1}{3} 
	\bigl( 4 \mathbf{v}_n - \mathbf{v}_{n-1} \bigr) + 
	\tfrac{2}{3} h \bigl( 2\mathbf{a}_n - \mathbf{a}_{n-1} 
	\bigr) + \mathcal{O} (h^4)
	\label{velocityPredictor} \\
	\intertext{followed with an evaluation of $\mathbf{a}^p_{n+1} = \mathbf{a} (t_{n+1}, \mathbf{x}^p_{n+1}, \mathbf{v}^p_{n+1})$.  The paired correctors belonging with these predictors are}
	\mathbf{x}_{n+1} & = \tfrac{1}{3} \bigl(
	4  \mathbf{x}_n - \mathbf{x}_{n-1} \bigr) +
	\tfrac{1}{24} h \bigl( \mathbf{v}^p_{n+1} +
	14 \mathbf{v}_n + \mathbf{v}_{n-1} \bigr) 
	\notag \\
	\mbox{} & \hspace{3.175cm} +
	\tfrac{1}{72} h^2 \bigl( 10 \mathbf{a}^p_{n+1} + 
	51 \mathbf{a}_n - \mathbf{a}_{n-1} \bigr) + 
	\mathcal{O} (h^4)
	\label{displacementCorrector} \\ 
	\mathbf{v}_{n+1} & = \tfrac{1}{3} 
	\bigl( 4 \mathbf{v}_n - \mathbf{v}_{n-1} \bigr) + 
	\tfrac{2}{3} h \mathbf{a}^p_{n+1} + \mathcal{O} (h^4)
	\label{velocityCorrector}
	\end{align}
\end{subequations}
which are followed with a re-evaluation for $\mathbf{a}_{n+1} = \mathbf{a} (t_{n+1}, \mathbf{x}_{n+1}, \mathbf{v}_{n+1})$.

Variables are to be updated according to $n \! - \! 1 \leftarrow n$ and $n \leftarrow n \! + \! 1$, after which counter $n$ gets incremented. Upon finishing the data management, a solution is ready for advancement to the next integration step, with looping continuing until $n = N$ whereat the solution becomes complete.

\section{Error and Step-Size Control}
\label{Sec:PI}

To be able to control the local truncation error one must first have an estimate for its value.  Here error is defined as a norm in the difference between predicted and corrected values.  A recipe for computing this is stated in the \textit{Strategy\/} section.  These expressions, although informative, cannot be used as stated because Taylor expansions for velocity have been applied to remove the next higher-order term in the Taylor series for displacement to improve accuracy.

An estimate for truncation error is simply
\begin{equation}
   \varepsilon_{n+1} = \frac{ \| \mathbf{x}_{n+1} - 
   	\mathbf{x}^p_{n+1} \| }
   {\max (1 , \| \mathbf{x}_{n+1} \| )}
   \label{truncationError}
\end{equation}
which can be used to control the size of a time step applied to an integrator, i.e., a local time step.  Our objective here is to keep $\varepsilon$ below some allowable error, i.e., a user specified tolerance denoted as \textit{tol}, typically set within the range of $[10^{-8}, 10^{-2}]$.

At this juncture it is instructive to introduce separate notations for the two time steps that arise in a typical implementation for an algorithm of this type into code.  Let $\Delta t$ denote the global time step, and let $h$ denote the local time step.  The global time step is considered to be uniformly sized at $\Delta t = T / N$, where $T$ is the time at which analysis stops and $N$ is the number of discrete nodes whereat information is to be passed back from the solver to its driver.  Typically $N$ is selected to be dense enough so that a user can create a suitable graphical representation of the result. On the other hand, the local time step $h$ that appears in  formul\ae\ (\ref{startUp1stOrderODEs}--\ref{pairedMethods}) is dynamically sized to maintain accuracy.  If error $\varepsilon$ becomes too large, then $h$ is reduced, and if it becomes too small, then $h$ is increased.  

If there is to be a local time step of size $h$ that adjusts dynamically, then the first question one must answer is: What is an acceptable value for $h$ to start an integration with?  It has been your professor's experience that the user is not as reliable in this regard as he\slash she would like to believe.  The following automated procedure has been found to be useful in this regard \cite{FreedIskovitz96}.  From the initial conditions, compute
\begin{displaymath}
   h_0 = \frac{ \| \mathbf{x}_0 \| }{\| \mathbf{v}_0 \|}
   \qquad \text{constrained so that} \qquad
   \frac{\Delta t}{100} < h_0 < \frac{\Delta t}{10}
\end{displaymath}
and with this initial estimate for the step size, take an Euler step forward $\mathbf{x}^p_1 = \mathbf{x}_0 + h_0 \mathbf{v}_0$, evaluate $\mathbf{v}^p_1 = \mathbf{v} (h_0, \mathbf{x}^p_1)$, follow with a trapezoidal correction $\mathbf{x}_1 = \mathbf{x}_0 + \tfrac{1}{2}  h_0 (\mathbf{v}^p_1 + \mathbf{v}_0)$, and re-evaluate  $\mathbf{v}_1 = \mathbf{v} (h_0, \mathbf{x}_1)$.  At this juncture, one can get an improved estimate for the initial step size via
\begin{displaymath}
	h_1 = 2 \left| \frac{\| \mathbf{x}_1 \| - \| \mathbf{x}_0 \|}
	{\| \mathbf{v}_1 \| + \| \mathbf{v}_0 \|} \right|
	\qquad \text{subject to} \qquad
   \frac{\Delta t}{1000} < h_1 .
\end{displaymath}
With this information, one can calculate the number of steps $S$ needed by a local solver to traverse the first step belonging to the global solver whose step size is $\Delta t$; specifically,
\begin{equation}
   S = \max \bigl( 2 , 
   \mathrm{round} ( \Delta t / h_1 ) \bigr) 
   \qquad \text{with} \qquad 
   h = \Delta t / S
   \label{initialStepSize}
\end{equation}
and a reasonable value for the initial, local, step size $h$ is now in hand.  As a minimum, there are to be two local steps taken for each global step traversed.

From here on a discrete PI controller (originally derived from control theory as a senior engineering project at Lund Institute of Technology in Lund, Sweden \cite{Gustafssonetal88}) is employed to automatically manage the size of $h$.  The goal of this PI controller is to allow a solution to traverse its path with maximum efficiency, all the while maintaining a specified tolerance on error.  

The P in PI stands for proportional feedback and accounts for current error, while the I in PI stands for integral feedback and accounts for an accumulation of error.  The simplest controller is an I controller.  For controlling step size, this I controller adjusts $h$ via \cite{Soderlind02}
\begin{displaymath}
C = \frac{h_{n+1}}{h_n} = \left( 
\frac{tol}{\varepsilon_{n+1}} \right)^{k_I}
\end{displaymath}
wherein $k_I$ designates gain in the integral feedback loop, while $tol$ is the maximum truncation error to be tolerated over a local step of integration.  Such controllers had been used by the numerical analysis community for a long time, and are known to be problematic \cite[pp.~31--35]{HairerWanner91}.  Controls engineers know that PI controllers are superior to I controllers, and for the task of managing $h$, in 1988 a team of students at Lund University derived \cite{Gustafssonetal88} 
\begin{displaymath}
C = \frac{h_{n+1}}{h_n} = \left( 
\frac{tol}{\varepsilon_{n+1}} \right)^{k_I+k_P} \left( 
\frac{\varepsilon_{n+1} \vphantom{l}}{tol} \right)^{k_P}
\end{displaymath}
wherein $k_P$ designates gain in the proportional feedback loop.  This PI controller has revolutionized how commercial-grade ODE solvers are built today.

A strategy for managing error by dynamically adjusting the size of time step $h$ can now be put forward.  To do so, it is instructive to introduce a second counter $s$ that decrements from $S$ down to 0.  It designates the number of steps left to go before reaching the node located at the end of a global step that the integrator is currently traversing.  $S$ needs to be redetermined each time the algorithm advances to its next global step.  If there is a discontinuity in step size $h$ across this interface, then the history variables will need to be adjusted using, e.g., a Hermite interpolator \cite{Shampine85}.  A suitable algorithm for controlling truncation error by managing step size is described below.

\medskip
Initialize the controller by setting $\varepsilon_n = 1$.
\begin{enumerate}
	\item After completing an integration for displacement $\mathbf{x}_{n+1}$, and possibly velocity $\mathbf{v}_{n+1}$, via any of the integrators given in Eqs.~(\ref{startUp1stOrderODEs}--\ref{pairedMethods}), calculate an estimate for its local truncation error $\varepsilon_{n+1}$ via Eq.~(\ref{truncationError}).
	\item Calculate a scaling factor $C$ that comes from the controller 
	\begin{displaymath}
	   C = \begin{cases}
	      \left( \vphantom{\frac{a}{a}} \right.
	      \frac{\mathit{tol}}{\varepsilon_{n+1}} 
	      \left. \vphantom{\frac{a}{a}} \right)^{0.7/(p+1)} 
	      \left( \frac{\varepsilon_n}{\mathit{tol}}
	      \right)^{0.4/(p+1)} &
	      \text{if} \; \varepsilon_{n} < \mathit{tol} 
	      \; \text{and} \; \varepsilon_{n+1} < \mathit{tol} \\
	      \left( \frac{\varepsilon_{n+1}}
	      {\mathit{tol}} \right)^{1/p} & \text{otherwise}
	   \end{cases}
	\end{displaymath}
	wherein \textit{tol\/} is the truncation error that the controller targets and $p$ is the order of the method, e.g., it appears as $\mathcal{O} (h^{p+1})$ in formul\ae\ (\ref{startUp1stOrderODEs}--\ref{pairedMethods}).
	\item If $C > 2$ plus $s > 3$ and $s$ is even, then double the step size $h = 2h$, halve the steps to go $s = s / 2$, and continue on to the next step.
	\item If $1 \leq C \leq 2$ then maintain the step size, decrement the counter, and on continue to the next step.
	\item If $C < 1$ yet $\varepsilon_{n+1} \leq \mathit{tol}$, then halve the step size $h = h / 2$, double the steps to go $s = 2 s$, and continue on to the next step.
	\item Else $C < 1$ and $\varepsilon_{n+1} > \mathit{tol}$, then halve the step size $h = h / 2$, double the steps to go $s = 2 s$, and \textit{repeat\/} the integration step from $n$ to $n \! + \! 1$.
\end{enumerate}
For the I controller, the gain on feedback has been set at $k_I = 1/p$.  For the PI controller, the gain on I feedback has been set at $k_I = 0.3 / (p+1)$ while the gain on P feedback has been set at $k_P = 0.4 / (p+1)$, wherein factors 0.3 and 0.4 have been selected based upon the developer's experience in working with their controller \cite{Gustafssonetal88,Soderlind02}.  By only admitting either a doubling or a halving of the current step size, a built-in mechanism is in play that mitigates the likelihood that wind-up or wind-down instabilities will happen in practice.

Whenever a step is to be halved, the displacement at a half step can be approximated via
\begin{equation}
\mathbf{x}_{n-\scriptfrac{1}{2}} = 
\tfrac{1}{2} ( \mathbf{x}_n + \mathbf{x}_{n-1} ) - 
\tfrac{1}{8} \, h ( \mathbf{v}_n - \mathbf{v}_{n-1} ) 
+ \mathcal{O} (h^4) + \mathcal{O} (h^{p+1})
\label{stepHalving} 
\end{equation}
which is a cubic Hermite interpolant \cite{Shampine85} whose accuracy is $\mathcal{O}(h^4)$ with $\mathcal{O}(h^{p+1})$ designating accuracy of the numerical method used to approximate displacements $\mathbf{x}_n$ and $\mathbf{x}_{n+1}$ and, for solvers (\ref{pairedStartUp} \& \ref{pairedMethods}), a like interpolation for velocities $\mathbf{v}_n$ and $\mathbf{v}_{n+1}$ will be required, too.

As a closing comment, many PECE methods are often implemented as $\text{PE}(\text{CE})^m$ methods with the correct\slash evaluate steps being repeated $m$ times, or until convergence.  It has been your professor's experience that PECE, i.e., $m=1$, is usually sufficient whenever the step size $h$ is properly controlled to keep the truncation error in check, provided a reasonable assignment for permissible error has been made, typically $tol \approx 10^{-(p+1)}$.

\section{Examples}

Examples are provided to illustrate the numerical methods put forward.  A chemical kinetics problem, popular in the numerical analysis literature \cite[pp.~115--116]{Haireretal93}, is considered for testing the two-step PECE method of Eqs.~(\ref{startUp1stOrderODEs} \& \ref{1stOrderODEs}) used to solve first-order systems of ODEs, including stiff ODEs.  The vibrational response of an formula SAE race car is simulated to illustrate a problem belonging to the class of solvers that are appropriate for applications of Newton's Second Law of motion.

\subsection{Brusselator}
\label{Sec:Brusselator}

The Brusselator describes a chemical kinetics problem where six substances are being mixed, and whose evolution through time is characterized by two, coupled, differential equations in two unknowns $A$ and $B$, viz.,
\begin{subequations}
    \begin{align}
    \dot{y}_1 & = A + y^2_1 y_2 - (B - 1) y_1 \notag \\
    \dot{y}_2 & = B y_1 - y^2_1 y_2 \notag
    \end{align}
\end{subequations}
whose eigenvalues are
\begin{displaymath}
    \lambda = \frac{1}{2} \left( - \left( 1 - B + A^2 \right) 
    \pm \sqrt{( 1 - B + A^2)^2 - 4 A^2} \right)
\end{displaymath}
where parameters $A$ and $B$ are, to an extent, at the disposal of a chemist. 

\begin{figure}
	{\par\centering
		\resizebox*{0.6\textwidth}{0.3\textheight}
		{\includegraphics{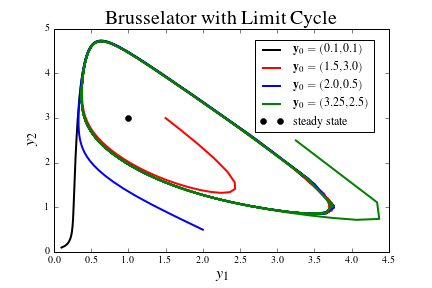}}
		\par}
	\caption{A concentration plot for a Brusselator response with $A=1$ and $B=3$.  Solutions are presented for several initial conditions.  All solutions approach a limit cycle.}
	\label{fig:brusselator1}
\end{figure}

This system exhibits vary different behaviors for different values of its parameters.  For values $A=1$ and $B=3$ (see Fig.~\ref{fig:brusselator1}) the solution converges to a limit cycle that orbits a steady-state attractor located at coordinate (1,~3) for these values of $A$ and $B$.  This limit cycle does not depend upon initial condition (IC), provided the IC does not reside at the steady state.  

\begin{figure}
	{\par\centering
		\resizebox*{0.9\textwidth}{0.225\textheight}
		{\includegraphics{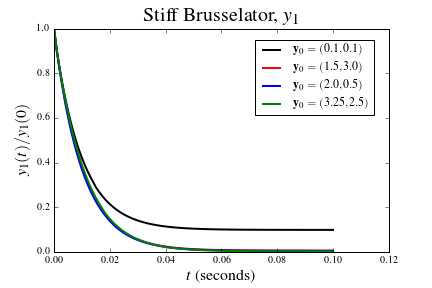}
		 \includegraphics{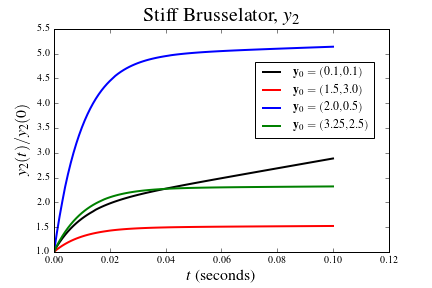}}
		\par}
	\caption{Brusselator response versus time with $A=100$ and $B=3$ for several initial conditions.  The response curves have been normalized against their initial values.}
	\label{fig:brusselator2}
\end{figure}

The behavior is very different for parameters $A=100$ and $B=3$.  Here the solutions rapidly settle in on asymptotic responses (see Fig.~\ref{fig:brusselator2}).  Figures \ref{fig:brusselator1} \& \ref{fig:brusselator2} came from the same system of equations, just different parameters.  

\begin{figure}
	{\par\centering
		\resizebox*{0.9\textwidth}{0.225\textheight}
		{\includegraphics{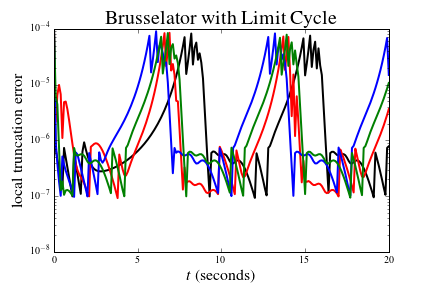}
		 \includegraphics{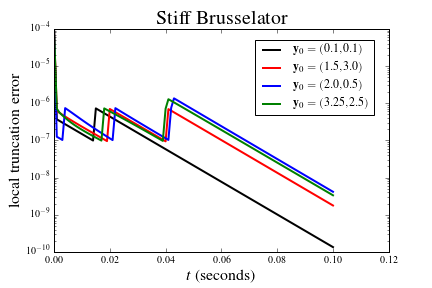}}
		\par}
	\caption{Local truncation error versus time for both Brusselator problems.  The error tolerance was set at $10^{-4}$, which is the upper horizontal axis in both plots.  The sawtooth response in the right plot was caused by the step size being doubled at those locations.  Oscillations of error in the left plot arose from the PI controller navigating corners.}
	\label{fig:brusselator3}
\end{figure}

Ability of the PI controller discussed in \S\ref{Sec:PI} to manage the local truncation error by adjusting the local step size is illustrated in Fig.~\ref{fig:brusselator3}.  Statistics gathered from these runs are reported on in Table~\ref{Table:Brusselator}

\begin{table}
	\small
	\begin{center}
    \begin{tabular}{|c|ccc|ccc|} \hline
    Initial &
    \multicolumn{3}{c|}{$A=1$, $B=3$, $t_{\mathrm{end}}=20$~s} &
    \multicolumn{3}{c|}{$A=100$, $B=3$, $t_{\mathrm{end}}=0.1$~s} \\
    \cline{2-7} Condition & \#steps & \#halved & \#doubled &
    \#steps & \#halved & \#doubled \\ \hline
    (0.1, 0.1) & 1186 & 6 & 9 & 353 & 0 & 6 \\
    (1.5, 3.0) & 1592 & 6 & 9 & 362 & 0 & 4 \\
    (2.0, 0.5) & 1332 & 7 & 10 & 467 & 0 & 3 \\
    (3.25, 2.5) & 1451 & 6 & 12 & 414 & 0 & 5 \\ \hline
    \end{tabular}
	\end{center}
	\normalsize
	\caption{Runtime statistics for the results plotted in Figs.~\ref{fig:brusselator1}--\ref{fig:brusselator3}.  There were 200 global steps for the limit cycle analyses, and 100 global steps for the stiff analyses.  In none of these numerical experiments did the integrator have to restart because of excessive error.}
	\label{Table:Brusselator}
\end{table}

The solutions in Fig.~\ref{fig:brusselator1} have an eigenvalue ratio of $| \lambda_{\max} | / | \lambda_{\min} | = 2.6$, whereas the solutions in Fig.~\ref{fig:brusselator2} have a ratio of $| \lambda_{\max} | / | \lambda_{\min} | = 9,602$.  Although there is no accepted `definition' for stiffness in the numerical analysis literature, there are some rules of thumb that exist.  Probably the simplest to apply is the ratio $\Lambda = | \lambda_{\max} | / | \lambda_{\min} |$ with $\Lambda \approx 10$ being the boundary.  Systems of ODEs whose ratio of extreme eigenvalues is less than about 10 do not exhibit stiffness; whereas, systems of ODEs whose ratio $\Lambda$ exceeds 10, and certainly 100, do exhibit stiffness.  

Explicit methods, e.g., the predictors presented herein, when used alone, do not fair well when attempting to acquire solutions from systems of ODEs that are mathematically stiff.  Implicit methods are needed, e.g., the correctors presented herein.  The solutions graphed in Fig.~\ref{fig:brusselator1} are for a non-stiff problem, while the solutions graphed in Fig.~\ref{fig:brusselator2} are for a stiff problem.  The implicit two-step method of Eqs.~(\ref{startUp1stOrderODEs} \& \ref{1stOrderODEs}) is a viable integrator for solving stiff systems of ODEs of first order; in contrast, explicit Runge-Kutta methods are not suitable.

\subsection{Vibrational Response of a Vehicle}

In this example we consider the vibrational response of a car as it travels down a roadway.  This response is excited by an unevenness in the roadway, accentuated by the speed of a vehicle.  This simulation determines the heave $z$, pitch $\theta$, and roll $\phi$ of a vehicle at its center of gravity excited by its traversal over a roadway.  

There are three degrees of freedom for this problem with the position $\mathbf{x}$, velocity $\mathbf{v}$, and acceleration $\mathbf{a}$ vectors taking on forms of
\begin{displaymath}
   \mathbf{x} = \left\{ \begin{matrix} 
      z \\ \theta \\ \phi
   \end{matrix} \right\} , \qquad
   \mathbf{v} = \left\{ \begin{matrix} 
   \dot{z} \\ \dot{\theta} \\ \dot{\phi}
   \end{matrix} \right\} , \qquad
   \mathbf{a} = \left\{ \begin{matrix} 
   \ddot{z} \\ \ddot{\theta} \\ \ddot{\phi}
   \end{matrix} \right\}
\end{displaymath}
wherein $\dot{z} = \partial{z} / \partial t$, $\ddot{z} = \partial^2 z / \partial t^2$, etc.  In our application of this simulator, we consider a formula SAE race car like the one our seniors design, fabricate and compete with every year in a cap stone project here at Texas~A\mbox{\&}M.

There are three matrices that establish the vibrational characteristics of a vehicle.  There is a mass matrix
\begin{equation*}
   \mathbf{M} = \begin{bmatrix}
	   m & 0 & 0 \\ 0 & J_{\theta} & 0 \\ 0 & 0 & J_{\phi}
   \end{bmatrix}
\end{equation*}
where $m$ is the collective mass of the car and its driver, $J_{\theta}$ is the moment of inertia resisting pitching motions, and $J_{\phi}$ is the moment of inertia resisting rolling motions. There is also a damping matrix
\begin{multline}
   \mathbf{C} = 
   \left[ \begin{matrix}
      c_1 + c_2 + c_3 + c_4 \\
      -(c_1 + c_2) \ell_f + (c_3 + c_4) \ell_r \\   
      -(c_1 - c_2) \rho_f + (c_3 - c_4) \rho_r  
   \end{matrix} \right. \notag \\  
   \left. \begin{matrix}
      -(c_1 + c_2) \ell_f + (c_3 + c_4) \ell_r &
         -(c_1 - c_2) \rho_f + (c_3 - c_4) \rho_r \\ 
      (c_1 + c_2) \ell_f^2 + (c_3 + c_4) \ell_r^2 &
         (c_1 - c_2) \ell_f \rho_f + (c_3 - c_4) \ell_r \rho_r  \\
      (c_1 - c_2) \ell_f \rho_f + (c_3 - c_4) \ell_r \rho_r &
         (c_1 + c_2) \rho_f^2 + (c_3 + c_4) \rho_r^2
   \end{matrix} \right] \notag
\end{multline}
and a like stiffness matrix
\begin{multline}
\mathbf{K} = 
\left[ \begin{matrix}
k_1 + k_2 + k_3 + k_4 \\
-(k_1 + k_2) \ell_f + (k_3 + k_4) \ell_r \\   
-(k_1 - k_2) \rho_f + (k_3 - k_4) \rho_r  
\end{matrix} \right. \notag \\  
\left. \begin{matrix}
-(k_1 + k_2) \ell_f + (k_3 + k_4) \ell_r &
-(k_1 - k_2) \rho_f + (k_3 - k_4) \rho_r \\ 
(k_1 + k_2) \ell_f^2 + (k_3 + k_4) \ell_r^2 &
(k_1 - k_2) \ell_f \rho_f + (k_3 - k_4) \ell_r \rho_r  \\
(k_1 - k_2) \ell_f \rho_f + (k_3 - k_4) \ell_r \rho_r &
(k_1 + k_2) \rho_f^2 + (k_3 + k_4) \rho_r^2
\end{matrix} \right] \notag
\end{multline}
wherein $c_1$ and $k_1$ are the effective damping coefficient and spring stiffness for the suspension located at the driver's front, $c_2$ and $k_2$ are located at the passenger's front, $c_3$ and $k_3$ are located at the passenger's rear, and $c_4$ and $k_4$ are located at the driver's rear.  Lengths $\ell_f$ and $\ell_r$ measure distance from the front and rear axles to the center of gravity (CG) for the car and driver with their sum being the wheelbase.  Lengths $\rho_f$ and $\rho_r$ measure distance from the centerline (CL) of the vehicle out to the center of a tire patch along the front and rear axles, respectively.  Typically, $\rho_f > \rho_r$ to allow a driver to take a tighter\slash shorter path into a corner during competition.

Interacting with these three matrices is a vector that establishes how a roadway excites a vehicle.  It is described by
\begin{displaymath}
   \mathbf{f} = \left\{ \begin{matrix}
   w - 
   c_1 \dot{R}_1 - c_2 \dot{R}_2 - c_3 \dot{R}_3 - c_4 \dot{R}_4
   - k_1 R_1 - k_2 R_2 - k_3 R_3 - k_4 R_4 \\
   \bigl( c_1 \dot{R}_1 + c_2 \dot{R}_2 
      + k_1 R_1 + k_2 R_2 \bigr) \ell_f - 
   \bigl( 3_1 \dot{R}_3 + c_4 \dot{R}_4 
      + k_3 R_3 + k_4 R_4 \bigr) \ell_r \\
   \bigl( c_1 \dot{R}_1 - c_2 \dot{R}_2 
      + k_1 R_1 - k_2 R_2 \bigr) \rho_f - 
      \bigl( 3_1 \dot{R}_3 - c_4 \dot{R}_4 
      + k_3 R_3 - k_4 R_4 \bigr) \rho_r
   \end{matrix} \right\}
\end{displaymath}
where $w$ is the weight (mass times gravity) of the car and its driver. Functions $R(t)$ and $\dot{R}(t)$ are for displacement and velocity occurring normal to a roadway, measured from smooth.  Roadway velocity is proportional to vehicle speed.  It is through these functions that time enters into a solution.  $R_i$ and $\dot{R_i}$, $i=1,2,3,4$, follow the same numbering scheme as the damping coefficients and spring stiffnesses.

To apply our numerical algorithm (\ref{pairedStartUp} \& \ref{pairedMethods}), one simply computes
\begin{displaymath}
   \mathbf{a} (t , \mathbf{x} , \mathbf{v}) = \mathbf{M}^{-1} 
   \cdot \bigl( \mathbf{f}(t) - \mathbf{C} \cdot \mathbf{v} 
   - \mathbf{K} \cdot \mathbf{x} \bigr)
\end{displaymath}
and assigns a suitable pair of ICs: one for displacement, and the other for velocity, as they pertain to the motion of a vehicle at its center of gravity.  Initial conditions can be cast is various ways.  The simplest ICs come from either starting at rest, or starting at a constant velocity on a smooth roadway.  Either way, one arrives at
\begin{displaymath}
   \mathbf{x}_0 = \mathbf{K}^{-1} \cdot \mathbf{f}_0 
   \quad \text{and} \quad 
   \mathbf{v}_0 = \left\{ \begin{matrix} 
      0 \\ 0 \\ 0
   \end{matrix} \right\}
   \quad \text{wherein} \quad
   \mathbf{f}_0 = \left\{ \begin{matrix} 
      w \\ 0 \\ 0
   \end{matrix} \right\}
\end{displaymath}
because $R_i = 0$ and $\dot{R}_i = 0$, $i=1,2,3,4$, in these two cases.  Remember, velocity $\mathbf{v}_0$ is not the speed of your car; rather, it is a change in vehicle motion with respect to its center of gravity.

\begin{figure}
	{\par\centering
		\resizebox*{0.9\textwidth}{0.225\textheight}
		{\includegraphics{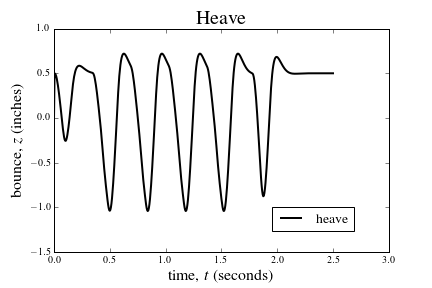}
			\includegraphics{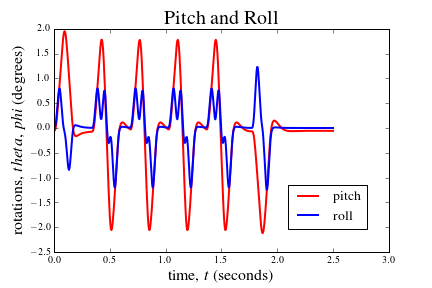}}
		\par}
	\caption{Heave $z$ is plotted against time in the left graphic, while pitch $\theta$ and roll $\phi$ are plotted against time in the right graphic.  Heave and pitch have static offsets, whereas roll does not.}
	\label{fig:fsae1}
\end{figure}

\begin{figure}
	{\par\centering
		\resizebox*{0.6\textwidth}{0.3\textheight}
		{\includegraphics{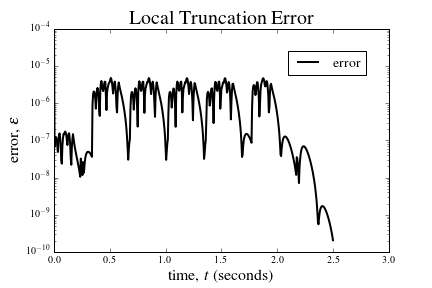}}
		\par}
	\caption{Local truncation error versus time for the FSAE race car driving over a sequence of bumps.  The error tolerance was set at $10^{-4}$, which is the upper horizontal axis of the plot. }
	\label{fig:fsae2}
\end{figure}

To illustrate the simulator, a roadway was constructed with five gradual waves at a wavelength equal to the wheelbase.  To excite roll, the passenger side lagged out of phase with the driver side by a tenth of the wheelbase.  Vehicle speed was set at 10~mph.  There were 500 global nodes so the density of output would produce nice graphs, for which there were 5,422 local integration steps required with 8 steps being doubled.  No steps were halved, and no steps required to be restarted.  The responses are plotted in Fig.~\ref{fig:fsae1}, while the errors are reported in Fig.~\ref{fig:fsae2}.  It is apparent that the integrator (\ref{pairedStartUp} \& \ref{pairedMethods}) performs to expectations, and that the PI controller of \S\ref{Sec:PI} does an admirable job in managing the local truncation error.

Vitals for the car that was simulated include: $m$ = 14~slugs ($w$ = 450~lbs), $J_{\theta} = 45 \text{ ft.lbs/(rad/sec}^2)$,  $J_{\phi} = 20 \text{ ft.lbs/(rad/sec}^2)$, $\ell_f$ = 3.2~ft, $\ell_r$ = 1.8~ft, $\rho_f$ = 2.1~ft, $\rho_r$ = 2~ft, the front dampers were set at 10 lbs/(in/sec) and the rears were set at 15 lbs/(in/sec), while the front springs had stiffnesses of 150 lbs/in and the rears were selected at 300~lbs/in.  These are reminiscent of a typical FSAE race car.

\section{Summary}

Two-step methods have been constructed that aspire to the structure of the well-known BDF2 formula.  A predictor is derived for each case allowing PECE solution schemes to be put forward.  The first method (\ref{startUp1stOrderODEs} \& \ref{1stOrderODEs}) that was introduced solves the classic problem where $\partial \mathbf{x} / \partial t = \mathbf{v} (t , \mathbf{x})$ subject to an IC of $\mathbf{x}(0) = \mathbf{x}_0$.  The second method (\ref{startup} \& \ref{PECE}) introduced solves a fairly atypical case where functions for both velocity $\mathbf{v} ( t , \mathbf{x})$ and acceleration $\mathbf{a} (t , \mathbf{x} , \mathbf{v} )$ are given and a solution for the displacement $\mathbf{x}$ is sought, subject to an initial condition $\mathbf{x}(0) = \mathbf{x}_0$.  And the third method (\ref{pairedStartUp} \& \ref{pairedMethods}) melds these two algorithms to construct a solver for the case where acceleration is given via a function $\mathbf{a} ( t , \mathbf{x} , \mathbf{v})$ from which solutions for both velocity $\mathbf{v}$ and displacement $\mathbf{x}$ are sought, subject to initial conditions of $\mathbf{x}(0) = \mathbf{x}_0$ and $\mathbf{v}(0, \mathbf{x}_0) = \mathbf{v}_0$.  A PI controller is used to manage the local truncation error by dynamically adjusting the size of the local time step.  All integrators have been illustrated using non-trivial example problems.

\newpage
\medskip\noindent\textbf{Acknowledgment}\medskip

The author is grateful to Prof.\ Kai Diethelm, Institut Computational Mathematics, Technische Universit\"at, Braunschweig, Germany for critiquing this document and for providing instructive comments.

\medskip\noindent\textbf{References}\medskip

%% References

\small
\bibliographystyle{spbasic}

%\bibliographystyle{asmems4}
%\bibliography{my}

%% \begin{thebibliography}{00}
%% \bibitem{}
%%
%% \end{thebibliography}
    
\end{document}